\documentclass[12pt]{article}
\usepackage{graphicx}
\usepackage[T1]{fontenc}
\usepackage{lmodern}
\usepackage{mathtools}
\usepackage{microtype}
\usepackage{xcolor}
\usepackage{enumitem}
\usepackage{hyperref}
\usepackage{amsmath,amsthm,amssymb,amsfonts,euscript,enumerate}
\usepackage{color,wrapfig}
\usepackage{pxfonts}
\usepackage{verbatim}
\numberwithin{equation}{section}
\newtheorem{theorem}{Theorem}[section]
\newtheorem{lemma}[theorem]{Lemma}
\newtheorem{proposition}[theorem]{Proposition}
\newtheorem{corollary}[theorem]{Corollary}
\theoremstyle{definition}

\theoremstyle{remark}

\newcommand{\R}{\mathbb R}

\newcommand{\dd}{\,\mathrm d}

\newcommand{\loc}{\mathrm{loc}}

\newcommand{\Cstar}{C_*}

\begin{document}
\title{Rigidity of singular eternal solutions of the fast diffusion equation}
\author{Shu-Yu Hsu\\
Department of Mathematics\\
National Chung Cheng University\\
168 University Road, Min-Hsiung\\
Chia-Yi 621, Taiwan, R.O.C.\\
e-mail: shuyu.sy@gmail.com}
\date{Nov 29, 2026}
\smallbreak \maketitle
\begin{abstract}
Let $n\ge 3$, $0<m<\frac{n-2}{n}$,  $\beta<0$ and $\alpha=\frac{2\beta}{1-m}$. For any $\lambda>0$, let  $U_{\lambda}(x,t)=e^{-\alpha t}f_{\lambda}(e^{-\beta t}x)$ be the singular eternal self-similar solution of 
\begin{equation*}
u_t=\frac{n-1}{m}\Delta u^m\quad\mbox{ in }(\mathbb{R}^n\setminus\{0\})\times(-\infty,\infty)
\end{equation*}
whose profile satisfies $\underset{\substack{r\to 0}}{\lim}\frac{r^2f_{\lambda}(r)^{1-m}}{\log r^{-1}}=\frac{2(n-1)(n-2-nm)}{|\beta|(1-m)}$ and $\underset{\substack{r\to\infty}}{\lim}r^{\frac{n-2}{m}}f_{\lambda}(r)=\lambda^{\frac{2}{1-m}-\frac{n-2}{m}}$. Let $u$ be a singular eternal solution of the fast diffusion equation. We prove that if the rescaled function $v(x,t)=e^{\alpha t}u(e^{\beta t}x,t)$
satisfies 
\begin{equation*}
0<v(x,t)\le f_{\bar\lambda}(x)\quad\forall (x,t)\in (\mathbb{R}^n\setminus\{0\})\times(-\infty,\infty)
\end{equation*}
for some constant $\bar\lambda>0$ and
\begin{equation*}
E_{\lambda_0}(t_j)=\int_{\mathbb{R}^n\setminus\{0\}} |v(x,t_j)-f_{\lambda_0}(x)||x|^{-\mu}\,dx
 \to 0\quad\mbox{ as }j\to\infty
\end{equation*}
for  some constant ${\lambda_0}>0$ and sequence $t_j\to-\infty$ as $j\to\infty$, then
\begin{equation*}
v(x,t)=f_{\lambda_0}(x)\quad\mbox{and}\quad
 u(x,t)=e^{-\alpha t}f_{\lambda_0}(e^{-\beta t}x)
 \quad\forall x\in \mathbb{R}^n\setminus\{0\}, t\in\mathbb{R}.
\end{equation*}
\end{abstract}

\vskip 0.2truein



\vskip 0.2truein
\setcounter{section}{0}

\section{Introduction}
\setcounter{equation}{0}
\setcounter{theorem}{0}

Let $n\ge 3$ and $0<m<\frac{n-2}{n}$. We will study the rigidity property of singular eternal solution of the equation,
\begin{equation}\label{eq:physicalpde}
 u_t=\frac{n-1}{m}\Delta u^m,\quad u>0,\quad\text{ in }(\mathbb{R}^n\setminus\{0\})\times\R.
\end{equation}
This equation was studied extensively 
by P. Daskalopoulos,  N. Sesum,  \cite{DS1}, \cite{DS2},  S.Y. Hsu \cite{Hs},  K.M. Hui \cite{H1}, \cite{H2}, etc. When  $n\ge 3$ and $m=\frac{n-2}{n+2}$, \eqref{eq:physicalpde} arises in the study of the Yamabe flow \cite{DS1}, \cite{DS2}. We refer the readers to the survey paper \cite{A} by D.G.~Aronson
and the books \cite{DK}, \cite{V}, by P. Daskalopoulos, C.E.~Kenig and J.L.~Vazquez on the recent results of \eqref{eq:physicalpde}.
 
Let
\begin{equation}\label{eq:parameters}
 n\geq3,\qquad 0<m<\frac{n-2}{n},\qquad
 \beta<0,\qquad \alpha=\frac{2\beta}{1-m},
\end{equation}
and
\begin{equation}\label{eq:exponents}
 a=\frac{2}{1-m},\qquad b=\frac{n-2}{m},\qquad
 \mu=n-\frac{2}{1-m},\qquad q=\frac1{1-m}.
\end{equation}
Then
\begin{equation}\label{eq:exponentfacts}
 2<a<n<b,\qquad 0<\mu<n-2,\qquad
 am=a-2,\qquad bm=n-2,\qquad b-a=\frac{\mu}{m}>\mu.
\end{equation}
For any $\lambda>0$, let $f_{\lambda}$ be the unique radially symmetric solution of 
\begin{equation}\label{eq:stationary}
\frac{n-1}{m}\Delta f_{\lambda}^m+\alpha f_{\lambda}+\beta x\cdot\nabla f_{\lambda}=0, \quad f_{\lambda}>0,\quad\mbox{ in }\mathbb{R}^n\setminus\{0\},
\end{equation}
given by Theorem 1.1 and Remark 1 of \cite{HP} which satisfies
\begin{equation}\label{eq:f-profile-asymptotic}
\lim_{r\downarrow0}\frac{r^2f_{\lambda}(r)^{1-m}}{\log r^{-1}}=\frac{2(n-1)(n-2-nm)}{|\beta|(1-m)}\quad\mbox{ and }\quad \lim_{r\to\infty}r^{\frac{n-2}{m}}f_{\lambda}(r)=\lambda^{\frac{2}{1-m}-\frac{n-2}{m}}.
\end{equation}
By \eqref{eq:f-profile-asymptotic} we have
\begin{align}
 f_\lambda(r)&\sim \Cstar^q r^{-a}\bigl(\log(1/r)\bigr)^q\quad\mbox{ as }r\downarrow0,\label{eq:profileorigin}\\
 f_\lambda(r)&\sim\lambda^{a-b}r^{-b}\qquad\qquad\mbox{ as }r\to\infty,\label{eq:profileinfty}
\end{align}
where
\begin{equation*}
 \Cstar=\frac{2(n-1)(n-2-nm)}{(1-m)|\beta|}>0.
\end{equation*}
In particular, for each fixed \(\lambda>0\), there exist
constants $C>0$, $R_0>0$ and $0<r_0<\min (R_0,e^{-1})$ such that
\begin{equation}\label{eq:profilebounds}
 f_\lambda(r)\leq
 \begin{cases}
 Cr^{-a}\bigl(\log(1/r)\bigr)^q,&0<r<r_0,\\
 Cr^{-b},&r>R_0.
 \end{cases}
\end{equation}

For any $\lambda>0$, let
\begin{equation}\label{eq:self-similar-formula}
U_{\lambda}(x,t)=e^{-\alpha t}f_{\lambda}(e^{-\beta t}x)\quad\forall x\in .
\end{equation}
Then $U_{\lambda}$ is an eternal self-similar solution of \eqref{eq:physicalpde} in $(\mathbb{R}^n\setminus\{0\})\times\mathbb{R}$. 
Suppose \(u\in C^{2,1}((\mathbb{R}^n\setminus\{0\})\times\R)\) is a solution of
\eqref{eq:physicalpde}. We define its rescaling as
\begin{equation}\label{eq:rescale}
 v(x,t)=e^{\alpha t}u(e^{\beta t}x,t)\quad\forall x\in \mathbb{R}^n\setminus\{0\},t\in \mathbb{R}.
\end{equation}
Then $v$ satisfies
\begin{equation}\label{eq:rescaledpde}
 v_t=\frac{n-1}{m}\Delta v^m+\alpha v+\beta x\cdot\nabla v\quad\mbox{ in }(\mathbb{R}^n\setminus\{0\})\times \mathbb{R}.
\end{equation}
For any \(\lambda>0\), let
\begin{equation}\label{eq:criticalE}
 E_\lambda(t)=\int_{\mathbb{R}^n\setminus\{0\}} |v(x,t)-f_\lambda(x)||x|^{-\mu}\dd x
 \quad\text{ in }[0,\infty].
\end{equation}
The definition allows an infinite value until finiteness is established.

We will prove the following results in this paper:

\begin{theorem}[Rigidity from a backward weighted limit]\label{thm:rigidity}
Suppose  \eqref{eq:parameters} holds. Let \(u\) be a solution of \eqref{eq:physicalpde} which satisfies
\begin{equation}\label{eq:upperclass}
 0<u(x,t)\leq U_{\bar\lambda}(x,t)
 \quad\forall x\in{\mathbb{R}^n\setminus\{0\}},\ t>0
\end{equation}  
for some constant \(\bar\lambda>0\).
Let $v$ be given by \eqref{eq:rescale}  and  $E_\lambda(t)$ be given by \eqref{eq:criticalE} for any $\lambda>0$ and $t\in\mathbb{R}$.
Suppose there exists a constant \(\lambda_0>0\) and a sequence \(t_j\to-\infty\) such that
\begin{equation}\label{eq:backwardweighted}
 E_{\lambda_0}(t_j)\to 0\quad\mbox{ as }j\to\infty.
\end{equation}
Then
\begin{equation}\label{eq:rigidityconclusion0}
 u(x,t)=e^{-\alpha t}f_{\lambda_0}(e^{-\beta t}x)
 \quad\forall (x,t)\in (\mathbb{R}^n\setminus\{0\})\times \mathbb{R}
\end{equation}
or equivalently
\begin{equation}\label{eq:rigidityconclusion}
v(x,t)=f_{\lambda_0}(x)\quad\forall (x,t)\in (\mathbb{R}^n\setminus\{0\})\times \mathbb{R}.
\end{equation}
\end{theorem}

Note that the condition \eqref{eq:upperclass} is equivalent to the condition,
\begin{equation}\label{eq:upperbound}
0<v(x,t)\leq f_{\bar\lambda}(x)
       \quad\mbox{ in }(\mathbb{R}^n\setminus\{0\})\times \mathbb{R}.
\end{equation}

\begin{theorem}[Rigidity from backward local convergence and tightness]
\label{thm:local}
Suppose  \eqref{eq:parameters} holds. Let \(u\) be a solution of \eqref{eq:physicalpde} which satisfies \eqref{eq:upperclass}
for some constant \(\bar\lambda>0\).  Let $v$ be given by \eqref{eq:rescale}.
Suppose there exists a constant \(\lambda_0>0\) and a sequence \(t_j\to-\infty\) such that
\begin{equation}\label{eq:localbackward}
 v(\cdot,t_j)\longrightarrow f_{\lambda_0}
 \quad\text{in }L^1_{\loc}(\mathbb{R}^n\setminus\{0\})\quad\mbox{ as }j\to\infty,
\end{equation}
and there exists a constant  \(t_*\in\R\) such that
\begin{equation}\label{eq:tightness}
\lim_{\rho\downarrow0}\ \sup_{t\leq t_*}
 \int_{0<|x|<\rho}|v(x,t)-f_{\lambda_0}(x)||x|^{-\mu}\dd x=0.
\end{equation}
Then \eqref{eq:rigidityconclusion0} or equivalently \eqref{eq:rigidityconclusion} holds. 
\end{theorem}

\begin{corollary}[An integrable-envelope criterion]\label{cor:envelope}
Suppose  \eqref{eq:parameters} holds. Let \(u\) be a solution of \eqref{eq:physicalpde} which satisfies \eqref{eq:upperclass}
for some constant \(\bar\lambda>0\).  Let $v$ be given by \eqref{eq:rescale}  and  $E_\lambda(t)$ be given by \eqref{eq:criticalE} for any $\lambda>0$ and $t\in\mathbb{R}$.
Suppose there exists a constant \(\lambda_0>0\) and a sequence \(t_j\to-\infty\) such that
\eqref{eq:localbackward} holds
and there exist a nonnegative measurable
function \(G\) on \(\{0<|x|<1\}\) and a constant  \(t_*\in\R\) such that

\begin{equation}\label{eq:envelope}
 \int_{0<|x|<1}G(x)|x|^{-\mu}\dd x<\infty,
 \qquad
 |v(x,t)-f_{\lambda_0}(x)|\leq G(x)
 \quad\forall 0<|x|<1,\ t\leq t_*.
\end{equation}
Then \eqref{eq:rigidityconclusion0} or equivalently \eqref{eq:rigidityconclusion} holds. 
\end{corollary}

We say that a function  \(u\) is a solution of \eqref{eq:physicalpde} if $u\in C^{2,1}((\mathbb{R}^n\setminus\{0\})\times\R)\) is positive on $(\mathbb{R}^n\setminus\{0\})\times\R$ and satisfies \eqref{eq:physicalpde} in the classical sense in $(\mathbb{R}^n\setminus\{0\})\times\R$. We will assume that $u$ is a solution of \eqref{eq:physicalpde} which satisfies \eqref{eq:upperclass} for some $\bar\lambda>0$ for the rest of the paper. We will also let $v$ be given by \eqref{eq:rescale} for the rest of the paper. Then \eqref{eq:upperbound} holds.
For any $a\in\mathbb{R}$, we let $\operatorname{sign}(a)=1$ if $a>0$,  $\operatorname{sign}(a)=-1$ if $a<0$ and  $\operatorname{sign}(a)=0$ if $a=0$.  For any  $x_0\in\R^n$ and $R>0$, we let $B_R(x_0)=\{x\in\mathbb{R}^d: |x-x_0|<R\}$ and $B_R=B_R(0)$. 

Since  \(\mu=n-\frac{2}{1-m}\),
\begin{equation}\label{eq:criticalcancel}
 \alpha-n\beta+\mu\beta=\beta\left(\frac{2}{1-m}-n+\mu\right)=0.
\end{equation}
Hence  by \eqref{eq:self-similar-formula}, \eqref{eq:rescale} and \eqref{eq:criticalE},
\begin{align} \label{eq:physicaldistance}
 \int_{\mathbb{R}^n\setminus\{0\}} |u(y,t)-U_{\lambda_0}(y,t)|\,|y|^{-\mu}\dd y
=&e^{-\alpha t}\int_{\mathbb{R}^n\setminus\{0\}}|v(e^{-\beta t}y,t)-f_{\lambda_0}(e^{-\beta t}y)||y|^{-\mu}\dd y\notag\\
 =&e^{-(\alpha-n\beta+\mu\beta)t}\int_{\mathbb{R}^n\setminus\{0\}}|v(x,t)-f_{\lambda_0}(x)||x|^{-\mu}\dd x\notag\\
=&E_{\lambda_0}(t).
\end{align}
This identity holds also as an identity of extended nonnegative
integrals. 

\section{Rigidity theorems for the eternal solutions}
\setcounter{equation}{0}
\setcounter{theorem}{0}

In this section we will prove Theorem \ref{thm:rigidity}, Theorem \ref{thm:local} and Corollary \ref{cor:envelope}. We first start with a contraction lemma.

Let $\lambda_0>0$ and 
\begin{equation}\label{eq:QB}
 h_{\lambda_0}=v-f_{\lambda_0},\qquad z_{\lambda_0}=v^m-f_{\lambda_0}^m,\qquad
 Q_{\lambda_0}=|h_{\lambda_0}|,\qquad B_{\lambda_0}=|z_{\lambda_0}|.
\end{equation}

\begin{lemma}[The difference inequality]\label{lem:kato}
The following inequality holds.
\begin{equation}\label{eq:kato}
 \partial_tQ_{\lambda_0}\leq\frac{n-1}{m}\Delta B_{\lambda_0}+\alpha Q_{\lambda_0}+\beta x\cdot\nabla Q_{\lambda_0}\quad\hbox{ in }\mathcal D'({\mathbb{R}^n\setminus\{0\}}\times(s,t))\quad\forall -\infty<s<t<\infty.
\end{equation}
\end{lemma}
\begin{proof}
This result was proved in Lemma 4.2 of \cite{HS} using Kato's inequality \cite{K}. For the sake of completeness we give an elementary proof of this lemma here.
By subtracting  \eqref{eq:stationary} from \eqref{eq:rescaledpde} we get
\begin{equation}\label{eq:h equation}
 h_{\lambda_0,t}=\frac{n-1}{m}\Delta z_{\lambda_0}+\alpha h_{\lambda_0}+\beta x\cdot\nabla h_{\lambda_0}\quad\hbox{ in }(\mathbb{R}^n\setminus\{0\})\times\mathbb{R}.
\end{equation}
Since \(r\mapsto r^m\) is strictly increasing on \((0,\infty)\),
\begin{equation}\label{eq:h-z-sign}
\operatorname{sign}h_{\lambda_0}=\operatorname{sign}z_{\lambda_0}
\end{equation}
Moreover
\begin{equation}\label{eq:ht-sign}
\operatorname{sign}(h_{\lambda_0})h_{\lambda_0}=|h_{\lambda_0}|,\quad\partial_t|h_{\lambda_0}|=\operatorname{sign}(h_{\lambda_0})h_{\lambda_0,t},
 \quad \nabla|h_{\lambda_0}|=\operatorname{sign}(h_{\lambda_0})\nabla h_{\lambda_0}
\end{equation}
for a.e. $(x,t)\in (\mathbb{R}^n\setminus\{0\})\times\mathbb{R}$. Let
\(\Phi_\varepsilon(r)=(r^2+\varepsilon^2)^{1/2}\). Then
\begin{equation}\label{eq:regularization of sign fcn}
 \Delta\Phi_\varepsilon(z_{\lambda_0})
 =\Phi_\varepsilon'(z_{\lambda_0})\Delta z_{\lambda_0}
  +\Phi_\varepsilon''(z_{\lambda_0})|\nabla z_{\lambda_0}|^2
 \geq\Phi_\varepsilon'(z_{\lambda_0})\Delta z_{\lambda_0}.
\end{equation}
Let $0\le \psi\in C_0^{\infty}((\mathbb{R}^n\setminus\{0\})\times\mathbb{R})$. 
Since \(z\) and its derivatives are bounded on the support of $\psi$, multiplying \eqref{eq:regularization of sign fcn} by $\psi$ and integrating over $(\mathbb{R}^n\setminus\{0\})\times\mathbb{R}$ we get
\begin{equation}\label{eq:regularized-ineqn}
\iint_{(\mathbb{R}^n\setminus\{0\})\times\mathbb{R}}\Phi_\varepsilon(z_{\lambda_0})\Delta\psi (z_{\lambda_0})\,dx
 \geq\iint_{(\mathbb{R}^n\setminus\{0\})\times\mathbb{R}}\Phi_\varepsilon'(z_{\lambda_0})\psi(z_{\lambda_0})\Delta z_{\lambda_0}\,dx
\end{equation}
Letting \(\varepsilon\downarrow0\) in \eqref{eq:regularized-ineqn},
\begin{equation*}
\iint_{(\mathbb{R}^n\setminus\{0\})\times\mathbb{R}}|z_{\lambda_0}|\Delta\psi(x)\,dx\geq\iint_{(\mathbb{R}^n\setminus\{0\})\times\mathbb{R}}\operatorname{sign}(z_{\lambda_0})\psi(x)\Delta z_{\lambda_0}\,dx.
\end{equation*}
Hence
\begin{equation}\label{eq;Delta-z-sign-ineqn}
\Delta|z_{\lambda_0}|\geq\operatorname{sign}(z_{\lambda_0})\Delta z_{\lambda_0}
 \quad\text{in distributions}.
\end{equation}
Multiplying \eqref{eq:h equation} by \(\psi\operatorname{sign}h_{\lambda_0}\) and integrating over $(\mathbb{R}^n\setminus\{0\})\times\mathbb{R}$,
by \eqref{eq:h-z-sign}, \eqref{eq:ht-sign} and \eqref{eq;Delta-z-sign-ineqn} we get
\eqref{eq:kato} and the lemma follows.
\end{proof}

For any \(0<\theta<\mu\), we let
\begin{equation}\label{eq:regularweight}
 E_{\lambda_0,\theta}(\tau)=\int_{(\mathbb{R}^n\setminus\{0\})\times \mathbb{R}}Q_{\lambda_0}(x,\tau)|x|^{-\mu+\theta}\dd x.
\end{equation}

\begin{lemma}[The regularized weighted estimate]\label{lem:theta}
For any \(0<\theta<\mu\) and $\lambda_0>0$, \(E_{\lambda_0,\theta}\) is finite and continuous
on $\mathbb{R}$ and
\begin{equation}\label{eq:thetaestimate}
E_{\lambda_0,\theta}(t)\leq e^{-\beta\theta(t-s)}E_{\lambda_0,\theta}(s)\quad\forall 
-\infty<s<t<\infty.
\end{equation}
\end{lemma}

\begin{proof}
By direct computation,
\begin{equation}\label{eq:lapweight}
 \Delta |x|^{-\mu+\theta}
 =-(\mu-\theta)(n-2-\mu+\theta)
       |x|^{-\mu+\theta-2}=-(\mu-\theta)\left(\frac{2m}{1-m}+\theta\right)
       |x|^{-\mu+\theta-2}<0,
\end{equation}
and
\begin{align}
 |\nabla |x|^{-\mu+\theta}|&=(\mu-\theta)|x|^{-\mu+\theta-1},\notag\\
 (\alpha-n\beta)|x|^{-\mu+\theta}-\beta x\cdot\nabla |x|^{-\mu+\theta}
 &=-\beta\theta |x|^{-\mu+\theta}.\label{eq:driftweight}
\end{align}
By \eqref{eq:profilebounds}, \eqref{eq:upperbound} and \eqref{eq:QB}, there exist constants $C>0$, $R_0>0$ and $0<r_0<\min (R_0,e^{-1})$ such that
\begin{equation}\label{eq:smallbounds}
 Q_{\lambda_0}(x,\tau)\leq Cr^{-a}\bigl(\log(1/r)\bigr)^q,
 \qquad B_{\lambda_0}(x,\tau)\leq Cr^{-am}\bigl(\log(1/r)\bigr)^{mq}
 \quad\forall 0<r<r_0, \tau\in\mathbb{R},
\end{equation}
and
\begin{equation}\label{eq:largebounds}
 Q_{\lambda_0}(x,\tau)\leq Cr^{-b},\qquad
 B_{\lambda_0}(x,\tau)\leq Cr^{-(n-2)}\quad\forall r>R_0,\tau\in\mathbb{R}.
\end{equation}
Hence by \eqref{eq:upperbound}, \eqref{eq:QB}, \eqref{eq:smallbounds} and \eqref{eq:largebounds},
\begin{equation*}
Q_{\lambda_0}(x,\tau)|x|^{n-1-\mu+\theta}\le \left\{\begin{aligned}
&C_1r^{n-1-a-\mu+\theta}\bigl(\log(1/r)\bigr)^q\quad\forall 0<r=|x|<r_0\\
&C_1\qquad\qquad\qquad\qquad\quad\forall r_0\le |x|\le R_0\\
&C_1r^{n-1-b-\mu+\theta}\qquad\qquad\quad\forall r=|x|>R_0
\end{aligned}\right.
\end{equation*}
for some constant $C_1>0$. Let
\begin{equation*}
H(r):=\left\{\begin{aligned}
&C_1r^{n-1-a-\mu+\theta}\bigl(\log(1/r)\bigr)^q\quad\forall 0<r<r_0\\
&C_1\qquad\qquad\qquad\qquad\quad\forall r_0\le r\le R_0\\
&C_1r^{n-1-b-\mu+\theta}\qquad\qquad\quad\forall r>R_0.
\end{aligned}\right.
\end{equation*}
Then
\begin{equation}\label{eq:Q-H-ineqn}
Q_{\lambda_0}(x,\tau)|x|^{n-1-\mu+\theta}\le H(r)\quad\forall r=|x|>0,\tau\in\mathbb{R}.
\end{equation}
Since
\[
\int_0^{r_0}H(r)\,dr=C_1\int_0^{r_0}r^{n-1-a-\mu+\theta}
             \bigl(\log(1/r)\bigr)^q\dd r
 =C\int_0^{r_0}r^{\theta-1}
             \bigl(\log(1/r)\bigr)^q\dd r<\infty
\]
and
\[\int_{R_0}^{\infty}H(r)\,dr=C_1\int_{R_0}^{\infty}r^{n-1-b-\mu+\theta}\dd r
 =C_1\int_{R_0}^{\infty}r^{a-b+\theta-1}\dd r<\infty,
\]
because \(\theta<\mu<b-a\), $H\in L^1(0,\infty)$. Then  by \eqref{eq:Q-H-ineqn} \(E_{\lambda_0,\theta}(t)<\infty\) for any $t\in\mathbb{R}$. By \eqref{eq:Q-H-ineqn} and the Lebesgue dominated convergence theorem  
\(E_{\lambda_0,\theta}(t)\) is a continuous  function of $t\in\mathbb{R}$.

Similar to the proof of Theorem 1.3 of \cite{HP} we now choose a smooth function \(\eta:[0,\infty)\to[0,1]\) equal to
\(0\) on \([0,1]\) and to \(1\) on \([2,\infty)\).
For any \(0<2\rho<R\), let
\[
 \chi_{\rho,R}(x)=\eta(|x|/\rho)\bigl(1-\eta(|x|/R)\bigr).
\]
Then $\chi_{\rho,R}(x)=0$ for \(|x|\leq\rho\) and \(|x|\geq2R\) and 
\(\chi_{\rho,R}(x)=1\) when \(2\rho\leq r\leq R\). Its derivatives are supported
on the two shells
\[
 A_\rho=\{\rho<r<2\rho\},\qquad A_R=\{R<r<2R\},
\]
where
\[
 |\nabla\chi_{\rho,R}|\leq C\ell^{-1},\qquad
 |\Delta\chi_{\rho,R}|\leq C\ell^{-2}
 \quad\text{on }A_\ell,\quad \ell=\rho,R.
\]
For brevity, put \(\chi=\chi_{\rho,R}\) and
\[
 F_{\rho,R}(\tau)=\int_{\mathbb{R}^n\setminus\{0\}} Q_{\lambda_0}(x,\tau)|x|^{-\mu+\theta}\chi(x)\dd x.
\]
Note that \(F_{\rho,R}\) is Lipschitz on each compact time interval.
Multiplying \eqref{eq:kato} by $|x|^{-\mu+\theta}\chi(x)$ and integrating by parts, by \eqref{eq:lapweight} and \eqref{eq:driftweight} we have
\begin{align}\label{eq:cutoffineq}
 &F_{\rho,R}'(\tau)\notag\\
 \leq&\frac{n-1}{m}\int_{\mathbb{R}^n\setminus\{0\}} B_{\lambda_0}\,\Delta(|x|^{-\mu+\theta}\chi)\dd x
   +\int_{\mathbb{R}^n\setminus\{0\}} Q_{\lambda_0}\bigl((\alpha-n\beta)|x|^{-\mu+\theta}\chi
                     -\beta x\cdot\nabla(|x|^{-\mu+\theta}\chi)\bigr)\dd x
   \notag\\
 \leq&-\beta\theta F_{\rho,R}(\tau)
   +\frac{n-1}{m}\int_{\mathbb{R}^n\setminus\{0\}} B_{\lambda_0}
         \bigl(|x|^{-\mu+\theta}\Delta\chi
                    +2\nabla |x|^{-\mu+\theta}\cdot\nabla\chi\bigr)\dd x\notag\\
&\quad   -\beta\int_{\mathbb{R}^n\setminus\{0\}} Q_{\lambda_0} |x|^{-\mu+\theta} x\cdot\nabla\chi\dd x\quad\mbox{ for a.e. }\tau\in\mathbb{R},0<2\rho< R\notag\\
\leq&-\beta\theta F_{\rho,R}(\tau)
   +\frac{n-1}{m}\int_{A_{\rho}} B_{\lambda_0}
         ||x|^{-\mu+\theta}\Delta\chi
                    +2\nabla |x|^{-\mu+\theta}\cdot\nabla\chi|\dd x\notag\\
&\quad+\frac{n-1}{m}\int_{A_R} B_{\lambda_0}
         ||x|^{-\mu+\theta}\Delta\chi
                    +2\nabla |x|^{-\mu+\theta}\cdot\nabla\chi|\dd x\notag\\
&\quad   -\beta\int_{A_{\rho}} Q_{\lambda_0} |x|^{-\mu+\theta} |x\cdot\nabla\chi|\dd x-\beta\int_{A_R} Q_{\lambda_0} |x|^{-\mu+\theta} |x\cdot\nabla\chi|\dd x\quad\forall 0<2\rho< R\mbox{ and  a.e. }\tau\in\mathbb{R}\notag\\
\leq&-\beta\theta F_{\rho,R}(\tau)+I_{\rho,R}. 
\end{align}
On each shell \(A_\ell\), $\ell=\rho, R$, \(|A_\ell|\le C\ell^n\) and
\[
 ||x|^{-\mu+\theta}\Delta\chi+2\nabla |x|^{-\mu+\theta}\cdot\nabla\chi|
 \leq C\ell^{-\mu+\theta-2},\qquad
 ||x|^{-\mu+\theta} x\cdot\nabla\chi|\leq C\ell^{-\mu+\theta}.
\]
Since \(n-\mu=a\) and \(a-am=2\), by \eqref{eq:smallbounds} we have
\begin{align}
 \int_{A_\rho}B_{\lambda_0}
  ||x|^{-\mu+\theta}\Delta\chi+2\nabla |x|^{-\mu+\theta}\cdot\nabla\chi|\dd x
 &\leq C\rho^{n-am-\mu+\theta-2}(\log(1/\rho))^{mq}\notag\\
 &=C\rho^\theta (\log(1/\rho))^{mq}\longrightarrow0\quad\mbox{ as }\rho\to 0
 \label{eq:innerdiffusion}
\end{align}
and
\begin{equation}\label{eq:innerdrift}
 \int_{A_\rho}Q_{\lambda_0}|x|^{-\mu+\theta}|x\cdot\nabla\chi|\dd x
 \leq C\rho^{n-a-\mu+\theta}(\log(1/\rho))^q
 =C\rho^\theta (\log(1/\rho))^q\longrightarrow0\quad\mbox{ as }\rho\to 0.
\end{equation}
By \eqref{eq:exponentfacts} and \eqref{eq:largebounds},
\begin{equation}\label{eq:outerdiffusion}
 \int_{A_R}B_{\lambda_0}
  ||x|^{-\mu+\theta}\Delta\chi+2\nabla |x|^{-\mu+\theta}\cdot\nabla\chi|\dd x
 \leq CR^{n-(n-2)-\mu+\theta-2}
 =CR^{-\mu+\theta}\longrightarrow0\quad\mbox{ as }R\to\infty,
\end{equation}
and
\begin{equation}\label{eq:outerdrift}
 \int_{A_R}Q_{\lambda_0}|x|^{-\mu+\theta}|x\cdot\nabla\chi|\dd x
 \leq CR^{n-b-\mu+\theta}
 =CR^{a-b+\theta}\le CR^{-\mu+\theta}\longrightarrow0\quad\mbox{ as }R\to\infty.
\end{equation}
By \eqref{eq:innerdiffusion}, \eqref{eq:innerdrift}, \eqref{eq:outerdiffusion} and \eqref{eq:outerdrift}, for any $0<2\rho<R$ there exists a constant $\varepsilon_{\rho,R}>0$
independent of $\tau$ such that
\begin{equation}\label{integral-upper-bd1}
I_{\rho,R}\le\varepsilon_{\rho,R}\to 0\quad\mbox{ as }\rho\to 0 \quad\mbox{ and }R\to\infty.
\end{equation}
Then by \eqref{eq:cutoffineq} and \eqref{integral-upper-bd1},
\begin{align}\label{eq:F-ineqn}
&F'_{\rho, R}(\tau)\le -\beta\theta F_{\rho, R}(\tau)+\varepsilon_{\rho,R}\quad\forall\tau\in\mathbb{R}\notag\\
\Rightarrow\quad&\frac{\partial}{\partial\tau}(e^{\beta\theta\tau}F_{\rho, R}(\tau))\le\varepsilon_{\rho,R}e^{\beta\theta\tau}\quad\forall\tau\in\mathbb{R}\notag\\
\Rightarrow\quad&F_{\rho, R}(t)\le e^{-\beta\theta(t-s)}F_{\rho, R}(s)+\frac{\varepsilon_{\rho,R}}{\beta\theta}\quad\forall -\infty<s<t<\infty.
\end{align}
Letting $\rho\to 0$ and $R\to\infty$ in \eqref{eq:F-ineqn}, by \eqref{integral-upper-bd1} we get \eqref{eq:thetaestimate} and the lemma follows.

\end{proof}

\begin{proposition}[Critical weighted contraction]\label{prop:critical}
If \(E_{\lambda_0}(s)<\infty\), then for any $t>s$,
\begin{equation}\label{eq:criticalcontraction}
E_{\lambda_0}(t)\leq E_{\lambda_0}(s).
\end{equation}
\end{proposition}

\begin{proof}
Let \(t>s\). We choose \(\theta_0\in(0,\mu)\) and let
\(0<\theta\leq\theta_0\). Then for any \(0<|x|\leq1\),
\[
 Q_{\lambda_0}(x,s)|x|^{-\mu+\theta}\leq Q_{\lambda_0}(x,s)|x|^{-\mu}.
\]
Since \(E_{\lambda_0}(s)<\infty\), the right-hand side is integrable in $B_1\setminus\{0\}$.
Hence by the Lebesgue dominated convergence theorem,
\begin{equation}\label{eq:E-limit1}
\lim_{\theta\downarrow0}\int_{0<|x|\le 1}Q_{\lambda_0}(x,s)|x|^{-\mu+\theta}\,dx=\int_{0<|x|\le 1}Q_{\lambda_0}(x,s)|x|^{-\mu}\,dx
\end{equation}
For \(|x|\geq1\),
\begin{equation}\label{eq:Q-estimate1}
 Q_{\lambda_0}(x,s)|x|^{-\mu+\theta}
 \leq Q_{\lambda_0}(x,s)|x|^{-\mu+\theta_0}.
\end{equation}
By the proof of Lemma~\ref{lem:theta}, the right hand side of \eqref{eq:Q-estimate1} is integrable in $\mathbb{R}^n\setminus B_1(0)$. Hence by the Lebesgue dominated convergence theorem,
\begin{equation}\label{eq:E-limit2}
\lim_{\theta\downarrow0}\int_{|x|> 1}Q_{\lambda_0}(x,s)|x|^{-\mu+\theta}\,dx=\int_{|x|>1}Q_{\lambda_0}(x,s)|x|^{-\mu}\,dx
\end{equation}
By \eqref{eq:E-limit1} and \eqref{eq:E-limit2},
\begin{equation}\label{eq:E-limit3}
 \lim_{\theta\downarrow0}E_{\lambda_0,\theta}(s)=E_{\lambda_0}(s).
\end{equation}
By Lemma \ref{lem:theta} \eqref{eq:thetaestimate} holds. Letting \(\theta\downarrow0\) in \eqref{eq:thetaestimate}, by \eqref{eq:E-limit3} and Fatou's lemma,
\[
 E_{\lambda_0}(t)
 \leq\liminf_{\theta\downarrow0}E_{\lambda_0,\theta}(t)
 \leq\lim_{\theta\downarrow0}
           e^{-\beta\theta(t-s)}E_{\lambda_0,\theta}(s)
 =E_{\lambda_0}(s)
\]
and \eqref{eq:criticalcontraction} follows.
\end{proof}

\begin{proof}[Proof of Theorem~\ref{thm:rigidity}]
Fix an arbitrary \(t\in\R\). Then there exists $j_0\in\mathbb{Z}^+$ such that 
\(t_j<t\) for any $j\ge j_0$. Then by Proposition~\ref{prop:critical},
\begin{equation}\label{eq:E-compare3}
 0\leq E_{\lambda_0}(t)\leq E_{\lambda_0}(t_j)\quad\forall j\ge j_0.
\end{equation}
Letting \(j\to\infty\) in \eqref{eq:E-compare3}, by \eqref{eq:backwardweighted}
we get 
\begin{align*}
&E_{\lambda_0}(t)=0\quad\forall t\in\mathbb{R}\notag\\
\Rightarrow\quad&v(x,t)\equiv f_{\lambda_0}(x)\quad\forall x\in\mathbb{R}^n\setminus\{0\}, t\in\mathbb{R}
\end{align*}
and \eqref{eq:rigidityconclusion} or equivalently \eqref{eq:rigidityconclusion0} follows.
\end{proof}

\begin{proof}[Proof of Theorem~\ref{thm:local}]
By \eqref{eq:profileinfty} and \eqref{eq:upperbound}, for sufficiently
large \(R\) and every \(t\in\R\),
\begin{align}
 \int_{|x|>R}|v(x,t)-f_{\lambda_0}(x)||x|^{-\mu}\dd x
 &\leq C\int_R^\infty r^{n-1-b-\mu}\dd r\notag\\
 &=C\int_R^\infty r^{a-b-1}\dd r
 \leq C_1R^{-(b-a)}.
 \label{eq:uniformoutertail}
\end{align}
Let \(\varepsilon>0\). Since \(b-a>0\), we can choose \(R>1\) sufficiently large such that the right-hand side of
\eqref{eq:uniformoutertail} is less than \(\varepsilon\).
By \eqref{eq:tightness} there exists $0<\rho<1$ such that  
\[
 \sup_{t\leq t_*}\int_{0<|x|<\rho}
      |v(x,t)-f_{\lambda_0}(x)||x|^{-\mu}\dd x<\varepsilon.
\]
By \eqref{eq:localbackward} there exists a constant $t_0<t_*$ such that
\[
 \int_{\rho\leq|x|\leq R}|v(x,t)-f_{\lambda_0}(x)||x|^{-\mu}\dd x<\varepsilon\quad\forall t\le t_0.
\]
Hence
\[
 E_{\lambda_0}(t)<3\varepsilon\quad\forall t\le t_0.
\]
Thus
\(E_{\lambda_0}(t)\to0\) as \(t\to-\infty\).
Theorem~\ref{thm:rigidity} now gives \eqref{eq:rigidityconclusion}.

\end{proof}

\begin{proof}[Proof of Corollary~\ref{cor:envelope}]
For \(0<\rho<1\), \eqref{eq:envelope} implies
\[
 \sup_{t\leq t_*}\int_{0<|x|<\rho}|v(x,t)-f_{\lambda_0}(x)||x|^{-\mu}\dd x
 \leq\int_{0<|x|<\rho}G(x)|x|^{-\mu}\dd x.
\]
The right-hand side tends to zero by integrability of \(G(x)|x|^{-\mu}\).
Thus \eqref{eq:tightness} holds, and Theorem~\ref{thm:local} applies. Thus the corollary follows.
\end{proof}

\end{document}